\documentclass[12pt]{article}
\usepackage{amssymb}
\usepackage{mathrsfs}
\makeatletter
 
 \@addtoreset{equation}{section}
\makeatother

\newcommand{\prf}{\noindent{\bf Proof. }}

\newcommand{\rem}{\noindent{\bf Remark. }}

\newcommand{\ff}[1]{%
    {\bf F}_{#1}}

\newcommand{\qed}{\hbox{\rule[-2pt]{5pt}{11pt}}}

\newtheorem{dfn}{Definition}[section]
\newtheorem{thm}[dfn]{Theorem}
\newtheorem{prop}[dfn]{Proposition}

\newtheorem{exam}[dfn]{Example}

\newtheorem{prob}[dfn]{Problem}

\begin{document}
%\begin{linenumbers}
%\begin{frontmatter}
\title{Extremal invariant polynomials not satisfying the Riemann hypothesis}
\author{Koji Chinen\footnotemark[1]}
\date{}
\maketitle

\begin{abstract}
Zeta functions for linear codes were defined by Iwan Duursma in 1999. They were generalized to the case of some invariant polynomials by the preset author. One of the most important problems is whether extremal weight enumerators satisfy the Riemann hypothesis. In this article, we show there exist extremal polynomials of the weight enumerator type which are invariant under the MacWilliams transform and do not satisfy the Riemann hypothesis. 
\end{abstract}%  words

\footnotetext[1]{Department of Mathematics, School of Science and Engineering, Kindai University. 3-4-1, Kowakae, Higashi-Osaka, 577-8502 Japan. E-mail: chinen@math.kindai.ac.jp}
\footnotetext{This work was supported by JSPS KAKENHI Grant Number JP26400028. }


%\footnotetext[1]{}

%%%%%%%%%%%%%%%%%%%%%%%%%%%%%%%%%%%%%%%%%%%%%%%%%%%%%%%%%%%%%
\section{Introduction}\label{section:intro}
Iwan Duursma \cite{Du1} discovered the zeta functions for linear codes and he developed their theory in \cite{Du2} -- \cite{Du4}. Later, the present author pointed out that we can define the zeta functions also for invariant polynomials not being related to linear codes and discussed some properties of them, including their Riemann hypothesis (see \cite{Ch1} and \cite{Ch2}). Undoubtedly, such invariant polynomials must be taken into account in order to investigate the zeta functions more in detail. With these preceding results in mind, we begin with the following definition: 
%%%%%%%%%%%%%%%%%%%%%%%%% Definition %%%%%%%%%%%%%%%%%%%%%%%%%
\begin{dfn}\label{dfn:zeta}
For any homogeneous polynomial of the form 
\begin{equation}\label{eq:homogen}
W(x,y)=x^n+\sum_{i=d}^n A_i x^{n-i}y^i \quad (A_i\in {\bf C},\ A_d\ne0)
\end{equation}
and $q\in{\bf R}$ ($q>0$, $q\ne 1$), there exists a unique polynomial $P(T)\in{\bf C}[T]$ of degree at most $n-d$ such that
\begin{equation}\label{eq:zeta_duursma}
\frac{P(T)}{(1-T)(1-qT)}(y(1-T)+xT)^n=\cdots +\frac{W(x,y)-x^n}{q-1}T^{n-d}+ \cdots.
\end{equation}
We call $P(T)$ and $Z(T)=P(T)/(1-T)(1-qT)$ the zeta polynomial and the zeta function of $W(x,y)$, respectively. 
\end{dfn}
One should note that we must assume $d, d^\perp\geq 2$ where $d^\perp$ is defined by 
$$W^{\sigma_q}(x,y)=\pm x^n + A_{d^\perp} x^{n-d^\perp} y^{d^\perp}+ \cdots,$$
when considering the zeta functions (see \cite[p.57]{Du2}). 
See \cite[Appendix A]{Ch2} for example, for an elementary proof of existence and uniqueness of $P(T)$. If $W(x,y)=W_C(x,y)$ is a weight enumerator of a linear code $C$, the number $q$ is the size of the finite field over which $C$ is defined. 

If $C$ is an $[n, n/2, d]$ self-dual code over the finite field ${\ff q}$, then the zeta polynomial $P(T)$ of the weight enumerator $W_C(x,y)$ satisfies the functional equation of the form 
\begin{equation}\label{eq:func_eq}
P(T)=P\left(\frac{1}{qT}\right)q^g T^{2g}\qquad (g=n/2+1-d).
\end{equation}
The equation (\ref{eq:func_eq}) comes from the fact that $W_C(x,y)$ is invariant under the MacWilliams transform 
\begin{equation}\label{eq:macwilliams}
\sigma_q=\frac{1}{\sqrt{q}}\left(\begin{array}{rr} 1 & q-1 \\ 1 & -1 \end{array}\right),
\end{equation}
where the action of a matrix $\sigma=\left(\begin{array}{cc} a & b \\ c & d \end{array}\right)$ on a polynomial $f(x,y)\in {\bf C}[x,y]$ is defined by $f^\sigma(x,y)=f(ax+by, cx+dy)$. From now on, we call a polynomial $f(x,y) \in {\bf C}[x,y]$ satisfying $f^\sigma(x,y)=f(x,y)$ a $\sigma$-{\it invariant polynomial}. We can show that (\ref{eq:func_eq}) is true for any $\sigma_q$-invariant  polynomial of the form (\ref{eq:homogen}), regardless of whether it is related to an existing code or not (see \cite[p.59]{Du2}). Thus it is appropriate to formulate the Riemann hypothesis for a $\sigma_q$-invariant polynomial in the following form ($q\in{\bf R}, q>0, q\ne 1$): 
%%%%%%%%%%%%%%%%%%%%%%%%% Definition %%%%%%%%%%%%%%%%%%%%%%%%%
\begin{dfn}\label{dfn:RH}
A $\sigma_q$-invariant polynomial $W(x,y)$ of the form (\ref{eq:homogen}) satisfies the Riemann hypothesis if all the roots of its zeta polynomial $P(T)$ lie on the circle $|T|=1/\sqrt{q}$. 
\end{dfn}
If $W(x,y)$ is the weight enumerator $W_C(x,y)$ of a self-dual code $C$, it is usually called ``the Riemann hypothesis for $C$''. There are examples of self-dual codes both satisfying and not satisfying the Riemann hypothesis (see \cite[Section 4]{Du3} for example). The following problem is well-known and important (\cite[Open Problem 4.2]{Du3}): 
%%%%%%%%%%%%%%%%%%%%%% problem %%%%%%%%%%%%%%%%%%%%%%
\begin{prob}[Duursma]\label{prob:duursma_extremal}
Prove or disprove that all extremal weight enumerators satisfy the Riemann hypothesis. 
\end{prob} 
An extremal weight enumerator is the weight enumerator of an ``extremal code'' which is a self-dual code attaining the equality in the Mallows-Sloane bound (see \cite[p.105]{Du4} or \cite[p.139]{Pl} for example). So far, we do not know an example of an extremal weight enumerator which does not satisfy the Riemann hypothesis. 

However, if we allow $q$ to be any positive real number other than one, removing the structure of the linear codes, then we find that there exist $\sigma_q$-invariant polynomials with some extremal property which do {\it not} satisfy the Riemann hypothesis. The purpose of this article is to show this fact by constructing explicit examples of such polynomials (see Section \ref{section:ext_non_rh}) and to give a solution in negative to Problem \ref{prob:duursma_extremal} in a broad sense. Such examples belong to the ring of $\sigma_q$-invariant polynomials for $q=6+2\sqrt{5}$. We can find the value $q$ by using the binomial moments. The algorithm is explained in Section \ref{section:remarks}. 

It should also be noted that there are examples of $\sigma_q$-invariant polynomials without an extremal property which satisfy the Riemann hypothesis. The author \cite{Ch2} provided plenty of such examples. 

These results suggest that the extremal property is neither a necessary nor a sufficient condition for the Riemann hypothesis. For similar results for the case of so-called formal weight enumerators, the reader is referred to \cite{Ch3}. 

The rest of the paper is organized as follows: in Section \ref{section:inv_ring}, we introduce a certain invariant polynomial ring and discuss some properties of them, including the extremal property. In Section \ref{section:ext_non_rh}, we show the existence of extremal $\sigma_q$-invariant polynomials not satisfying the Riemann hypothesis in the ring which is introduced in Section \ref{section:inv_ring}. In the last section, we give some remarks and problems. 

In what follows, $[x]$ means the greatest integer not exceeding $x$ for a real number $x$. The Pochhammer symbol $(a)_n$ means $(a)_n=a(a+1)\cdots (a+n-1)$ for $n\geq 1$ and $(a)_0=1$. We put 
\begin{equation}\label{eq:basic_inv_q}
W_{2,q}(x,y)=x^2+(q-1)y^2
\end{equation}
for $q>0$, $q\ne 1$. 
%%%%%%%%%%%%%%%%%%%%%%%%%%%%%%%%%%%%%%%%%%%%%%%%%%%%
%%%%%%%%%%%%%%%%%%%%%%%%%%%%%%%%%%%%%%%%%%%%%%%%%%%%
\section{Construction of a certain invariant polynomial ring}\label{section:inv_ring}
Let $q=6+2\sqrt{5}$ and 
\begin{eqnarray*}
\sigma &=& \sigma_q=\frac{1}{\sqrt{q}}\left(\begin{array}{rr} 1 & q-1 \\ 1 & -1 \end{array}\right)\\
&=& \frac{1}{1+\sqrt{5}}\left(\begin{array}{cc} 1 & 5+2\sqrt{5} \\ 1 & -1 \end{array}\right).
\end{eqnarray*}
We also put
$$\tau=\left(\begin{array}{rr} 1 & 0 \\ 0 & -1 \end{array}\right).$$
We form a group $G_{6+2\sqrt{5}}$ in $GL_2({\bf C})$ generated by $\sigma$ and $\tau$ and construct a polynomial ring 
$$R_{6+2\sqrt{5}}={\bf C}[x,y]^{G_{6+2\sqrt{5}}}$$
which is invariant under the action of $G_{6+2\sqrt{5}}$. 
%%%%%%%%%%%%%%%%%%%%%% proposition %%%%%%%%%%%%%%%%%%%%%%
\begin{prop}\label{prop:group}
Let $G_{6+2\sqrt{5}}=\langle \sigma, \tau \rangle$. Then we have the following: \\
(i) $|G_{6+2\sqrt{5}}|=10$. \\
(ii) The Molien series of $G_{6+2\sqrt{5}}$ are
$$\Phi(\lambda)=\frac{1}{(1-\lambda^2)(1-\lambda^5)}.$$
\end{prop} 
\prf (i) First we note that $\sigma^2=\tau^2=I=\left(\begin{array}{rr} 1 & 0 \\ 0 & 1 \end{array}\right)$. We can also verify that $\tau\sigma$ has order five and so $\sigma\tau=(\tau\sigma)^4$. Thus we can see that 
$$G_{6+2\sqrt{5}}=\{(\tau\sigma)^i \tau^j\ ; \ 0\leq i \leq 4, j=0,1 \}
=\langle \tau\sigma \rangle \rtimes \langle \tau \rangle.$$
We have $|G_{6+2\sqrt{5}}|=10$ from this decomposition. 

\medskip
\noindent (ii) We can prove this by a direct calculation of the Molien series (see \cite[p.600]{MaSl}) 
$$\Phi(\lambda)=\frac{1}{|G_{6+2\sqrt{5}}|}\sum_{A\in G_{6+2\sqrt{5}}} \frac{1}{\det(I-\lambda A)}.$$
Indeed, the contribution of $A=\tau\sigma$ and $(\tau\sigma)^4$ to the sum is $4/(2\lambda^2+(1-\sqrt{5})\lambda+2)$, that of $A=(\tau\sigma)^2$ and $(\tau\sigma)^3$ is $4/(2\lambda^2+(1+\sqrt{5})\lambda+2)$, that of $A=(\tau\sigma)^5=I$ is $1/(1-\lambda)^2$ and other elements ($(\tau\sigma)^i \tau$, $0\leq i \leq 4$) contribute by $5/(1-\lambda^2)$. We get the desired formula by summing them up. \qed

\medskip
\noindent Proposition \ref{prop:group} (ii) tells us that the ring $R_{6+2\sqrt{5}}$ has two generators, one of which is of degree two and another is of degree five (see \cite[p.601]{MaSl} for a similar discussion). In fact, we can verify that 
\begin{eqnarray}
W_{2, 6+2\sqrt{5}}(x,y) &=& x^2+(5+2\sqrt{5})y^2, \label{eq:gen_deg2}\\
\psi_5(x,y) &=& x^5-(50+20\sqrt{5})x^3y^2+(225+100\sqrt{5})xy^4 \label{eq:gen_deg5}
\end{eqnarray}
have the properties
\begin{eqnarray*}
W_{2, 6+2\sqrt{5}}(x,y)^\sigma &=& W_{2, 6+2\sqrt{5}}(x,y)^\tau = W_{2, 6+2\sqrt{5}}(x,y),\\
\psi_5(x,y)^\sigma &=& \psi_5(x,y)^\tau = \psi_5(x,y).
\end{eqnarray*}
Thus we have proved the following theorem: 
%%%%%%%%%%%%%%%%%%%%%% theorem %%%%%%%%%%%%%%%%%%%%%%
\begin{thm}\label{thm:generators}
$$R_{6+2\sqrt{5}}={\bf C}[W_{2, 6+2\sqrt{5}}(x,y), \psi_5(x,y)].$$
\end{thm} 
\rem (i) The members $W(x,y)$ of $R_{6+2\sqrt{5}}$ have the property ``divisible by two'', that is, if we assume $W(x,y)$ is of the form (\ref{eq:homogen}), we have ``$A_i \ne 0$ $\Rightarrow$ $2|i$''. This is due to the invariance under $\tau$ and is obvious from (\ref{eq:gen_deg2}) and (\ref{eq:gen_deg5}). 

\medskip
\noindent(ii) It is not always true that there is a $\sigma_q$-invariant polynomial divisible by two other than powers of $W_{2,q}(x,y)$. A method for finding a suitable $q$ is explained in Section \ref{section:remarks}. 

\medskip
The polynomials  $W(x,y)\in R_{6+2\sqrt{5}}$ of the form (\ref{eq:homogen}) are 
\begin{equation}\label{eq:form_of_W}
W_{2, 6+2\sqrt{5}}(x,y)^l \psi_5(x,y)^m\qquad(l,m\geq 0,\quad (l,m)\ne (0,0))
\end{equation}
and their suitable linear combinations. If there are several polynomials of the form (\ref{eq:form_of_W}) of the same degree, we can cancel the terms of small powers of $y$ and can get a $W(x,y)$ with a large $d$ by taking a suitable linear combination. 
%%%%%%%%%%%%%%%%%%%%%% exam %%%%%%%%%%%%%%%%%%%%%%
\begin{exam}\label{exam:deg10}\rm 
In the case where $\deg W(x,y)=10$, possible values of $l,m$ in (\ref{eq:form_of_W}) are $(l,m)=(5,0),(0,2)$. So we can form $W(x,y)$ from two polynomials $W_{2, 6+2\sqrt{5}}(x,y)^5$ and $\psi_5(x,y)^2$. We have
\begin{eqnarray}
W_{10}^{\rm E}(x,y) &:=& \frac{4}{5}W_{2, 6+2\sqrt{5}}(x,y)^5 + \frac{1}{5}\psi_5(x,y)^2\nonumber\\
 &=& x^{10}+(1350+600\sqrt{5})x^6y^4-(5100+2280\sqrt{5})x^4y^6\nonumber\\
 & & +(36225+16200\sqrt{5})x^2y^8+(30500+13640\sqrt{5})y^{10}.\label{eq:deg10extremal}
\end{eqnarray}
Thus we can eliminate the term of $y^2$ and can get a polynomial with $d=4$. Obviously, there is no other polynomials with $d\geq 4$ at this degree, so it should be called the extremal polynomial of degree 10 (see Definition \ref{dfn:def_extremal} later). 
\end{exam} 
In the case of even degrees, we can show that it does not happen that more coefficients accidentally cancel. This is an analog of the Mallows-Sloane bound: 
%%%%%%%%%%%%%%%%%%%%%% theorem %%%%%%%%%%%%%%%%%%%%%%
\begin{thm}\label{thm:analog_mallows-sloane}
Suppose that $W(x,y)\in R_{6+2\sqrt{5}}$ is of the form (\ref{eq:homogen}) and $\deg W(x,y)$ is even. Then the following inequality holds: 
\begin{equation}
d\leq 2\left[\frac{n}{10}\right]+2.
\end{equation}
\end{thm} 
\prf We apply the method of \cite[p.624-628]{MaSl}. So we use a similar notation and state an outline only. Let 
\begin{eqnarray*}
W_2(x,y) &=& W_{2, 6+2\sqrt{5}}(x,y),\\
W'_{10}(x,y) &=& (\psi_5(x,y)^2-W_{2, 6+2\sqrt{5}}(x,y)^5)/(25(5+2\sqrt{5}))\\
 &=& y^2(x^2-y^2)^2 (x^2-(9+4\sqrt{5})y^2)^2.
\end{eqnarray*}
Then $W(x,y)\in R_{6+2\sqrt{5}}$ of the form (\ref{eq:homogen}) and of even degree is written in the form
\begin{equation}\label{eq:lin_comb}
W(x,y) = \sum_{r=0}^\mu a_r W_2(x,y)^{5\mu+\nu-5r} W'_{10}(x,y)^r \quad (\mu \geq 0),
\end{equation}
here we have $n=\deg W(x,y)=2(5\mu + \nu)$. Suppose we choose suitable $a_r$ and we cancel as many coefficients as possible. The right hand side of (\ref{eq:lin_comb}) is a linear combination of $\mu + 1$ polynomials, so we can at least make $y^2, y^4, \cdots, y^{2\mu}$ disapear. So we assume
\begin{equation}\label{eq:lin_comb_cancel}
W(x,y)=x^n + \sum_{r=\mu+1}^{5\mu+\nu} A_{2r} x^{n-2r} y^{2r}.
\end{equation}
Now we prove $A_{2\mu +2}\ne 0$. We replace $x$ by 1 and $y^2$ by $x$ in $W_2(x,y)$ and $W'_{10}(x,y)$. We put
\begin{eqnarray*}
f(x) &=& 1+(5+2\sqrt{5})x,\\
g(x) &=& x(1-x)^2 (1-(9+4\sqrt{5})x)^2.
\end{eqnarray*}
The function ${\mit\Phi}(x)=xf(x)^5/g(x)$ satisfies the conditions of the B\"urmann-Lagrange Theorem (see \cite[Chapter 19, Theorem 14]{MaSl}) and we can conclude that 
\begin{eqnarray}
A_{2\mu+2} &=& - \frac{1}{(\mu+1)!} \frac{d^\mu}{dx^\mu} 
               \left.\left\{\left( \frac{d}{dx} \frac{1}{f(x)^{5\mu+\nu}}\right) {\mit\Phi}(x)^{\mu+1}\right\}\right|_{x=0}\nonumber\\
           &=& \frac{(5\mu+\nu)(5+2\sqrt{5})}{(\mu+1)!(9+4\sqrt{5})^{2\mu+2}} \frac{d^\mu}{dx^\mu} 
              \left.\left\{ \frac{(1+(5+2\sqrt{5})x)^{4-\nu}}{(x-1)^{2\mu+2}(x-1/(9+4\sqrt{5}))^{2\mu+2}} \right\}\right|_{x=0}.\label{eq:coeff_A_2m+2}
\end{eqnarray}
Let 
$$F_\mu (x)=(x-\alpha)^{-2\mu-2}(x-\beta)^{-2\mu-2}$$
for $\alpha, \beta>0$ and $\mu \in {\bf Z}$, $\mu\geq 0$. Then it is easy to see that 
$$F^{(l)}_{\mu}(0)=\sum_{r=0}^l {{l}\choose{r}} (2\mu+2)_{l-r}(2\mu+2)_r \alpha^{-2\mu-2-l+r} \beta^{-2\mu-2-r}>0$$
for all $l\geq 0$. Moreover, since $4-\nu\geq 0$,  $\{(1+(5+2\sqrt{5})x)^{4-\nu}\}^{(l)}|_{x=0}>0$ unless $\{(1+(5+2\sqrt{5})x)^{4-\nu}\}^{(l)}$ is identically zero. Thus we can see that $A_{2\mu+2}>0$ for all $\mu\geq 0$ and that $d\leq 2\mu+2$. Using $n=2(5\mu+\nu)$ and the fact that $d$ is even, we get the desired inequality. \qed
%%%%%%%%%%%%%%%%%%%%%% exam %%%%%%%%%%%%%%%%%%%%%%
\begin{exam}\label{exam:A_4}\rm 
Let $\mu=1$ and $\nu=0$. Then Theorem \ref{thm:analog_mallows-sloane} gives the coefficient $A_{2\mu+2}=A_4$ of $W_{10}^{\rm E}(x,y)$ in Example \ref{exam:deg10}. It is indeed
\begin{eqnarray*}
A_4 &=& \frac{5(5+2\sqrt{5})}{2(9+4\sqrt{5})^4}\frac{d}{dx}
\left.\left\{\frac{(1+(5+2\sqrt{5})x)^4}{(x-1)^4(x-1/(9+4\sqrt{5}))^4}\right\}\right|_{x=0}\\
 &=& 1350+600\sqrt{5}
\end{eqnarray*}
and coincides with the relevant coefficient in (\ref{eq:deg10extremal}). 
\end{exam} 

\medskip
\noindent Now we define the extremal property in the case of even degrees: 
%%%%%%%%%%%%%%%%%%%%%% definition %%%%%%%%%%%%%%%%%%%%%%
\begin{dfn}\label{dfn:def_extremal} 
A polynomial $W(x,y)\in R_{6+2\sqrt{5}}$ of the form (\ref{eq:homogen}) and of even degree is called extremal if 
\begin{equation}\label{eq:def_extremal}
d = 2\left[\frac{n}{10}\right]+2
\end{equation}
is satisfied. 
\end{dfn} 
\rem In the case of odd degrees, there seems to be a similar bound to (\ref{eq:def_extremal}). The conjectural bound is 
\begin{equation}\label{eq:mallows-sloane_odd}
d\leq 2\left[\frac{n-5}{10}\right]+2
\end{equation}
(the proof seems to be difficult). 
%%%%%%%%%%%%%%%%%%%%%%%%%%%%%%%%%%%%%%%%%%%%%%%%%%%%
%%%%%%%%%%%%%%%%%%%%%%%%%%%%%%%%%%%%%%%%%%%%%%%%%%%%
\section{Extremal invariant polynomials not satisfying the Riemann hypothesis}\label{section:ext_non_rh}
We first prove the following: 
%%%%%%%%%%%%%%%%%%%%%% theorem %%%%%%%%%%%%%%%%%%%%%%
\begin{thm}\label{thm:RH_extremal_deg10}
The extremal $\sigma_q$-invariant polynomial $W_{10}^{\rm E}(x,y)$ of (\ref{eq:deg10extremal}) does not satisfy the Riemann hypothesis. 
\end{thm} 
\prf We can calculate the zeta polynomial $P_{10}^{\rm E}(T)$ of $W_{10}^{\rm E}(x,y)$ (by some computer algebra system) as follows: 
\begin{eqnarray}
P_{10}^{\rm E}(T) &=& \frac{1}{7}\left\{ (520+232\sqrt{5})T^4 - (320+144\sqrt{5})T^3 \right. \nonumber\\
& & \left. - (168+76\sqrt{5})T^2 - (30+14\sqrt{5})T + 5+2\sqrt{5} \right\}. \label{eq:zeta_deg10_extremal}
\end{eqnarray}
It is not hard to see that $P_{10}^{\rm E}(99/100)=-122341/2500000-(40154771/87500000)\sqrt{5}<0$ and $P_{10}^{\rm E}(1)=1>0$ (the latter formula holds for zeta polynomials in general, see \cite[p.59, formula (7)]{Du2}). It follows that $P_{10}^{\rm E}(T)$ has a real root $\alpha$ with $99/100 < \alpha < 1$, which does not lie on the circle $|T|=1/\sqrt{6+2\sqrt{5}}\approx 0.3090$. \qed

\medskip
In the case of odd degrees, the bound (\ref{eq:mallows-sloane_odd}) is not established yet, but once the degree is specified, we can determine the extremal $\sigma_q$-invariant polynomial in $R_{6+2\sqrt{5}}$ (we call $W(x,y)$ of the form (\ref{eq:homogen}) extremal if $d$ is the largest at the given degree $n$). The simplest case is degree five, where $\psi_5(x,y)$ in (\ref{eq:gen_deg5}) is itself extremal. Indeed, there are no other $W(x,y)$ of the form (\ref{eq:homogen}) with $d\geq 2$ at this degree. 
%%%%%%%%%%%%%%%%%% theorem %%%%%%%%%%%%%%%%%%
\begin{thm}\label{thm:6+2sqrt5}
The extremal $\sigma_q$-invariant polynomial $\psi_5(x,y)$ does not satisfy the Riemann hypothesis. 
\end{thm}
\prf The zeta polynomial $P_5(T)$ of $\psi_5(x,y)$ is given by
$$P_5(T)=-\frac{2+\sqrt{5}}{4}(4T-1+\sqrt{5})( 8T^2-4\sqrt{5}T +3-\sqrt{5} ).$$
The factor $8T^2-4\sqrt{5}T +3-\sqrt{5}$ has two distinct real roots which are not equal to $\pm 1/\sqrt{6+2\sqrt{5}}$. \qed
%%%%%%%%%%%%%%%%%%%%%%%%%%%%%%%%%%%%%%%%%%%%%%%%%%%%
%%%%%%%%%%%%%%%%%%%%%%%%%%%%%%%%%%%%%%%%%%%%%%%%%%%%
\section{Some remarks and problems}\label{section:remarks}
(I) {\it How to find suitable values of $q$}

\medskip
As was mentioned before, it is not always true that there is a $\sigma_q$-invariant polynomial divisible by two other than powers of $W_{2,q}(x,y)$ (it is well known and is easy to show that $W_{2,q}(x,y)$ is $\sigma_q$-invariant for any $q$ ($q>0, q\ne 1$), see \cite{Sl} for example). One efficient method to find values like $q=6+2\sqrt{5}$ is to use the binomial moments (see MacWilliams-Sloane \cite[p.131, Problem (6)]{MaSl}): 
%%%%%%%%%%%%%%%%%% theorem %%%%%%%%%%%%%%%%%%
\begin{thm}[Binomial moments]\label{thm:moment} Let 
\begin{equation}\label{eq:w(x,y)general}
W(x,y)=\sum_{i=0}^n A_i x^{n-i}y^i \qquad(A_0=1)
\end{equation}
be a $\sigma_q$-invariant polynomial. Then we have 
\begin{equation}
\sum_{i=0}^{n-\nu} {{n-i}\choose{\nu}} A_i = q^{\frac{n}{2}-\nu} \sum_{i=0}^\nu {{n-i}\choose{n-\nu}} A_i
\qquad(\nu=0,1,\cdots, n).
\end{equation}
\end{thm}
We consider the case where $\deg W(x,y)$ in (\ref{eq:w(x,y)general}) is odd and is divisible by two (see Remark (i) after Theorem \ref{thm:generators}). Then, since $A_{2i+1}=0$ ($0\leq i \leq (n+1)/2$), we can assume
\begin{equation}\label{eq:we_odd_deg}
W(x,y)=\sum_{i=0}^n A_i x^{2(n-i)+1} y^{2i}\qquad(A_0=1)
\end{equation}
by omitting $A_{2i+1}$ and renumbering the suffices of the coefficients $A_i$. Note that $\deg W(x,y)=2n+1$. Then we can easily rewrite (\ref{thm:moment}) in the form 
\begin{equation}\label{eq:linear_eq_we_div2}
\sum_{i=0}^{[\frac{2n+1-\nu}{2}]} {{2n+1-2i}\choose{\nu}} A_i - q^{n-\nu+1/2} \sum_{i=0}^{[\frac{\nu}{2}]} {{2n+1-2i}\choose{2n+1-\nu}} A_i =0
\qquad(\nu=0,1,\cdots, 2n+1). 
\end{equation}
The formula (\ref{eq:linear_eq_we_div2}) gives $2n+2$ linear equations of $A_0, A_1, \cdots, A_n$, but the cases $\nu=n+1, n+2, \cdots, 2n+1$ are essentially the same as the cases $\nu=n, n-1, \cdots, 0$, respectively. So it suffices to consider the cases $\nu=0,1, \cdots, n$. Thus we get a system of $n+1$ homogeneous linear equations of $n+1$ unknowns $A_0, A_1, \cdots, A_n$. It is necessary that this system of linear equations has a non-trivial solution since $A_0=1\ne 0$ for the existence of a $\sigma_q$-invariant polynomial $W(x,y)$ of degree $2n+1$. Let $C(n,q)$ be the coefficient matrix of the first $n+1$ equations (the cases $\nu=0,1, \cdots, n$) given by (\ref{eq:linear_eq_we_div2}). We can find candidates of $q$ from the solutions of $|C(n,q)|=0$. Let us consider the case $n=2$. Then we have
$$C(2,q)=
\left(\begin{array}{ccc} 1-q^2\sqrt{q} & 1 & 1 \\
                        5(1-q\sqrt{q}) & 3 & 1 \\
                        10(1-\sqrt{q}) & 3-\sqrt{q} & 0\end{array}\right)$$
and $|C(2,q)|=-(t-1)^3(t+2)(t^2-2t-4)$ ($\sqrt{q}=t$). We get $t=1+\sqrt{5}$ from $|C(2,q)|=0$, and it gives $q=(1+\sqrt{5})^2=6+2\sqrt{5}$. This value was obtained by the consideration above. 

\medskip
\rem (i) If $n=1$, 
$$C(1,q)=\left(\begin{array}{cc} 1-q\sqrt{q} & 1 \\ 3(1-\sqrt{q}) & 1 \end{array}\right)$$
and $|C(1,q)|=-(\sqrt{q}+2)(\sqrt{q}-1)^2$, so we get no non-trivial invariant polynomial. 

\medskip
\noindent (ii) We can also apply the same argument in the case of formal weight enumerators divisible by two (see \cite{Ch3}). 

\bigskip
\noindent (II) {\it Some problems on the ring $R_{6+2\sqrt{5}}$}

\medskip
It is desirable to establish an analog of the Mallows-Sloane bound for the case of odd degrees: 
%%%%%%%%%%%%%%%%%% problem %%%%%%%%%%%%%%%%%%
\begin{prob}\label{prob:mallows-sloane_odd}
Can one prove 
$$d\leq 2\left[\frac{n-5}{10}\right]+2$$
(see (\ref{eq:mallows-sloane_odd})) for $W(x,y)\in R_{6+2\sqrt{5}}$ of the form (\ref{eq:homogen}) and of odd degree ?
\end{prob}
The expression of $A_{2\mu+2}$ is too complicated to show $A_{2\mu+2}\ne0$ in this case. Other than the method of \cite[p.624-628]{MaSl}, Duursma \cite{Du4} developed the theory of invariant differential operators on invariant polynomial rings to give an alternative proof of the Mallows-Sloane bound. It is uncertain whether we can find suitable differential operators for the ring $R_{6+2\sqrt{5}}$. 

The second problem is concerned to the Riemann hypothesis: 
%%%%%%%%%%%%%%%%%% problem %%%%%%%%%%%%%%%%%%
\begin{prob}\label{prob:RH_6+2sqrt5}
Is there a $W(x,y)\in R_{6+2\sqrt{5}}$ of the form (\ref{eq:homogen}) which satisfies the Riemann hypothesis ?
\end{prob}
The author has not found an example of a $W(x,y)\in R_{6+2\sqrt{5}}$ satisfying the Riemann hypothesis. It is a striking difference from the rings of weight enumerators of self-dual codes, which contain plenty of $W(x,y)$ satisfying it. The ring $R_{6+2\sqrt{5}}$ is similar to them, in the sense that there is an inequality like the Mallows-Sloane bound, but the situation about the Riemann hypothesis is quite different. It would be tempting to ask the following question: 
%%%%%%%%%%%%%%%%%% problem %%%%%%%%%%%%%%%%%%
\begin{prob}\label{prob:equivalent}
Find an equivalent condition for a $\sigma_q$-invariant polynomial to satisfy the Riemann hypothesis. 
\end{prob}

\bigskip
\noindent{\it Acknowledgement. }This work was established mainly during the author's stay at University of Strasbourg for the overseas research program of Kindai University. He would like to express his sincere gratitude to  Professor Yann Bugeaud at University of Strasbourg for his hospitality and to Kindai University for giving him a chance of the program.

%\end{linenumbers}
\end{document}